\documentclass{ifacconf}

\usepackage{amsmath}
\usepackage{amssymb}
\usepackage{amsfonts}
\usepackage{latexsym}
\usepackage{graphicx}      
\usepackage{natbib}        

\graphicspath{ {./} }

\newcommand{\N}{{\mathbb N}}

\newcommand{\cV}{{\cal V}}

\begin{document}
\begin{frontmatter}

\title{Finding Attracting Sets using Combinatorial Multivector Fields\thanksref{footnoteinfo}} 

\thanks[footnoteinfo]{The research of T.W.\ was partially supported by the Simons Foundation
under Award~581334.}

\author[First]{Justin Thorpe} 
\author[Second]{Thomas Wanner} 

\address[First]{Department of Mathematical Sciences, George Mason University,
    Fairfax, Virginia 22030, USA (e-mail: jthorpe3@gmu.edu)}
\address[Second]{Department of Mathematical Sciences, George Mason University,
Fairfax, Virginia 22030, USA (e-mail: twanner@gmu.edu)}

\begin{abstract}                
We discuss the identification of attracting sets using combinatorial multivector fields (CMVF) from
Conley-Morse-Forman theory. A CMVF is a dynamical system induced by the action of a continuous
dynamical system on a phase space discretization that can be represented as a Lefschetz complex.
There is a rich theory under development establishing the connections between the induced and
underlying dynamics and emphasizing computability. We introduce the main ideas behind this theory
and demonstrate how it can be used to identify regions of interest within the global dynamics
via graph-based algorithms and the connection matrix.
\end{abstract}


\end{frontmatter}

%
%
\section{Introduction}

Early work by \cite{eidenschink:95a} and \cite{boczko:02a,boczko:etal:07a} showed how the geometry of a
triangulation of the phase space of a dynamical system could be leveraged to approximate
isolated invariant sets. Classical Conley theory could then be deployed to make statements about the interior dynamics.
The primary requirement on the decomposition has been that boundaries of the subsets must be non-parallel
to local flow, thus ensuring isolation of the sets of interest.
This requirement has remained a hurdle for algorithmic implementation.

More recent work by \cite{mrozek:17a,lipinski:etal:23a} and others has abstracted much of the
classical theory to the finite, discrete setting of Conley-Morse-Forman theory.
The foundation of this theory rests on the notion of a combinatorial multivector field (CMVF), which serves
as a means to represent the discrete dynamics of the decomposition.

Our paper covers recent work that aims to relax the above-mentioned flow transversality requirement
of the earlier approaches. We employ the more general setting of a Lefschetz complex to represent
the decomposition, and we gain the ability to sidestep the transversality issue by "absorbing" any
problematic facets into the CMVF in a theoretically consistent fashion.

We begin by introducing core concepts related to Lefschetz complexes and their homology, 
followed by the CMVF and its dynamics.
We briefly cover Morse decompositions in this context as well as the connection matrix,
which serves as a useful algebraic tool to organize the dynamical interactions within the decomposition.
Finally, we demonstrate our approach by extracting attracting sets for examples in
two and three dimensions.

%
%
\section{Lefschetz Complexes}

\begin{defn}
   Let~$F$ denote an arbitrary field. Then a pair~$(X,\kappa)$ is
   called a {\em Lefschetz complex} over~$F$ if $X$ is a finite set
   with $\N_0$-gradation $X = (X_k)_{k \in \N_0}$, and the mapping
   $\kappa : X \times X \to F$ is such that
   \begin{displaymath}
      \kappa(x,y) \neq 0
      \quad\mbox{ implies }\quad
      x\in X_k \;\;\mbox{ and }\;\; y \in X_{k-1} \; ,
   \end{displaymath}
   and such that for any $x,y \in X$ one has
   \begin{equation} \label{kappaprop}
      \sum_{z \in X} \kappa(x,z) \kappa(z,y) = 0 \; .
   \end{equation}
\end{defn}

The elements $x\in X$ are referred to as {\em cells\/}, the scalar
value~$\kappa(x,y)$ is called the {\em incidence coefficient\/} of
the cells~$x$ and~$y$, and the map~$\kappa$ is the {\em incidence
coefficient map\/}. In addition, one defines the {\em dimension\/}
of a cell $x\in X_k$ as~$k$, and denotes it by~$k = \dim x$. Whenever
the incidence coefficient map is clear from context, we often just
refer to~$X$ as the Lefschetz complex.

The incidence coefficient map induces a {\em boundary map\/} $\partial$,
which is defined on cells via the relation
\begin{equation} \label{defbnd}
   \partial x = \sum_{y \in X} \kappa(x,y) y
\end{equation}
This map sends a cell $x$ of dimension $k$ to a specific 
linear combination of cells of dimension $k-1$, called the
{\em boundary\/} of $x$.
The boundary map can be extended to map a linear combination
of $k$-dimensional cells to the corresponding linear combination
of the separate boundaries. One can leverage this extended form
to rewrite the summation condition in the definition of a Lefschetz
complex in the equivalent form
\begin{equation} \label{bndprop}
   \partial( \partial x) = 0 \quad \text{ for all cells }\quad x \in X
\end{equation}
\noindent
i.e., the boundary of any cell is itself boundaryless.

There are several other concepts that will be important for our
discussion of Lefschetz complexes.
A {\em facet\/} of a cell $x \in X$ is any cell $y$ which
satisfies $\kappa(x,y) \neq 0$. 
One can define a partial order on the cells of $X$ by letting
$x \le y$ if and only if for some integer $n \in \mathbb{N}$
there exist cells $x = x_1, \ldots, x_n = y$ such that $x_k$
is a facet of $x_{k+1}$ for all $k = 1, \ldots, n-1$.
This partial order on $X$ is called the {\em face relation\/}. 
Thus, if $x \le y$ then $x$ is called a {\em face\/} of $y$.

A subset $C \subset X$ of a Lefschetz complex is called {\em closed\/},
if for every $x \in C$ all the faces of the cell $x$ are also
contained in the subset $C$.
The {\em closure\/} of a subset $C \subset X$, denoted by $\mathrm{cl}\, C$,
is the collection of all faces of all cells in $C$.
Thus, a subset of a Lefschetz complex is closed if and only if it
equals its closure.
A subset $S \subset X$ is called {\em locally closed\/}, if its {\em mouth\/}
$\mathrm{mo}\, S = \mathrm{cl}\, S \setminus S$ is closed. Note
that every closed set is automatically locally closed, but the
reverse implication is generally false.

The following results from \cite{mrozek:batko:09a}, \cite{mrozek:wanner:25a},
and \cite{lipinski:etal:23a}, respectively, characterize the relationship
between a {\em Lefschetz subcomplex\/} and a {\em locally closed\/} subset.
\begin{thm}
    Let $(X,\kappa)$ be a given Lefschetz complex over a field~$F$.
    Then $S \subset X$ is again a Lefschetz complex,
    with respect to the restriction of $\kappa$ to $S \times S$,
    if the subset $S$ is locally closed.
\end{thm}
\begin{cor}
  A subset $S \subset X$ is locally closed, if and only if it is the
  difference of two closed subsets of $X$.
\end{cor}
\begin{cor}
  A subset $S \subset X$ is locally closed, if and only if it is an
  interval with respect to the face relation on $X$, i.e., whenever
  we have three cells with $S \ni x \le y \le z \in S$, then one
  has to have $y \in S$ as well.
\end{cor}

Lefschetz complexes are a very general mathematical concept encompassing
other types of complex.
For example, a {\em simplicial complex\/}, {\em cubical complexes\/}, and
any {\em regular CW complex\/} are all types of Lefschetz complex.
See \cite{massey:91a} and \cite{dlotko:etal:11a} for more details.

%
%
\section{Lefschetz Homology}

For both the definition of the Conley index, and for our discussion
of connection matrices, we need to briefly recall the concept of homology.
Lefschetz homology generalizes both {\em simplicial homology\/} as 
described in \cite{munkres:84a}, and {\em cubical homology\/} in the sense
of \cite{kaczynski:etal:04a}. In order to fix our notation, we provide
a brief introduction in the following. For more details, see \cite{lefschetz:42a}.

\begin{defn}
   The {\em $k$-th chain group\/} $C_k(X)$ consists of all formal linear
   combinations of $k$-dimensional cells with coefficients in the underlying field $F$.
\end{defn}

\noindent
The collection of all chain groups is $C(X) = (C_k(X))_{k \in {\mathbb Z}}$, where we let
$C_k(X) = \{ 0 \}$ for all $k < 0$ and $k > \dim X$.

\begin{defn}
   The {\em $k$-th boundary map\/} $\partial_k$ is a linear map from
   $C_k(X)$ to $C_{k-1}(X)$ which acts, by using~\eqref{defbnd},
   on chain groups as follows:
   \begin{displaymath}
      \partial \left( \sum_{i=1}^m \alpha_i \sigma_i \right) =
      \sum_{i=1}^m \alpha_i \partial\sigma_i
      \;\;\text{ for }\;\;
      \sum_{i=1}^m \alpha_i \sigma_i \in C_k(X)
   \end{displaymath}
\end{defn}

We can use $\partial$ to form a sequence of vector spaces and linear
maps between them in the following way.

\begin{defn}
   A {\em chain complex\/} $(C(X), \partial)$ is a sequence of chain groups $C_k(X)$
   and boundary maps $\partial_k$ in the form
   \begin{displaymath}
      \ldots C_{k+1}(X) \stackrel{\partial_{k+1}}{\longrightarrow}
      C_{k}(X) \stackrel{\partial_{k}}{\longrightarrow}
      \ldots \stackrel{\partial_{1}}{\longrightarrow}
      C_0(X) \stackrel{\partial_{0}}{\longrightarrow} \{ 0 \}
      \ldots   
   \end{displaymath}
   such that $\partial_{k} \circ \partial_{k+1} = 0$ for all $k \in \mathbb{Z}$.
\end{defn}

Recall that the last statement regarding the composition of two successive 
boundary operators in the context of an underlying Lefschetz complex
follows from the properties of~$\kappa$, see also~\eqref{kappaprop}
and~\eqref{bndprop}. We highlight two important subspaces of the $k$-th
chain group $C_k(X)$:

\begin{itemize}
   \item The elements of the subspace $Z_k(X) = \mathrm{ker}\;
   \partial_k$ are called the {\em $k$-cycles\/} of $X$.
   \item The elements of $B_k(X) = \mathrm{im}\;
   \partial_{k+1}$ are called the {\em $k$-boundaries\/} of $X$,
   and one has $B_k(X) \subset Z_k(X)$.
\end{itemize}

\begin{defn}
   The {\em $k$-th homology group\/} of the Lefschetz complex $X$
   is the quotient space 
   \begin{displaymath}
      H_k(X) \; = \;
      Z_k(X) / B_k(X) \; = \;
      \mathrm{ker}\;\partial_k / \mathrm{im}\;\partial_{k+1}
   \end{displaymath}
\end{defn}

\noindent
which is well-defined in view of $B_k(X) \subset Z_k(X)$.

Note that $H_k(X)$ is again a vector space over $F$.
The dimension of $H_k(X)$ is called the {\em $k$-th Betti number of~$X$\/}
and abbreviated as $\beta_k(X) = \dim H_k(X)$. One can show that the $k$-th
Betti number counts the number of independent $k$-dimensional holes in~$X$.

Finally, we need to consider the homology of a complex~$X$ relative to
a subcomplex $X_0 \subset X$. In this case, we consider the relative chain
groups $C_k(X,X_0) = C_k(X)/C_k(X_0)$. One can then show that~$\partial_k$
induces a relative boundary map $\partial_k : C_k(X,X_0) \rightarrow
C_{k-1}(X,X_0)$ via $\partial_k [x] = [\partial_k x]$, which turns the
relative chain groups into a chain complex.

\begin{defn} {\em Relative Homology\/}
   \begin{itemize}
      \item The elements of the subspace $Z_k(X,X_0) = \mathrm{ker}\;
      \partial_k$ are called the {\em relative $k$-cycles\/} of $(X,X_0)$.
      \item The elements of the subspace $B_k(X,X_0) = \mathrm{im}\;
      \partial_{k+1}$ are called the {\em relative $k$-boundaries\/} of $(X,X_0)$.
      \item The vector space $H_k(X,X_0) = Z_k(X,X_0)/B_k(X,X_0)$
      is called the {\em $k$-th relative homology group\/} of $(X,X_0)$.
   \end{itemize}
\end{defn}

One can think of relative homology as a measure of holes in the complex~$X$,
after~$X_0$ has been contracted to a point.

%
%
\section{Multivector Fields via Phase Space Discretizations}

A {\em combinatorial multivector field\/} derives from, and extends,
the concept of {\em combinatorial vector field\/} introduced
in~\cite{forman:98a, forman:98b}. Multivector fields were developed
by \cite{mrozek:17a} and \cite{lipinski:etal:23a}, providing
more flexibility for encoding dynamical behavior. Going forward,
we assume that~$X$ is a Lefschetz complex over a field~$F$; typically
the rational numbers~$\mathbb{Q}$ or a finite field.
Then a {\em multivector field\/} on~$X$ can be defined as follows:
\begin{defn}
   A {\em multivector field\/} $\mathcal{V}$ on a Lefschetz complex
   $X$ is a partition of $X$ into locally closed sets.
\end{defn}

Recall that a subset $V \subset X$ is called locally closed if its mouth
$\mathrm{mo}\, V = \mathrm{cl}\, V \setminus V$ is closed, where
closedness in turn is defined via the face relation in the Lefschetz
complex. This implies that for each multivector $V \in \mathcal{V}$
the relative homology $H_*(\mathrm{cl}\, V, \mathrm{mo}\, V)$
is well-defined, and it allows for the following classification
of multivectors:
\begin{itemize}
   \item A multivector is {\em critical\/} if
   $H_*(\mathrm{cl}\, V, \mathrm{mo}\, V) \neq 0$.
   \item A multivector is {\em regular\/} if
   $H_*(\mathrm{cl}\, V, \mathrm{mo}\, V) = 0$.
\end{itemize}
Since a multivector is locally closed, it is itself
a Lefschetz complex, and in fact one has
$H_*(V) \cong H_*(\mathrm{cl}\, V, \mathrm{mo}\, V)$. At first
glance the distinction between regular and critical multivectors
might seem strange. Note, however, that in classical Conley theory
the notion of isolating block can be used to identify isolated
invariant sets. Loosely speaking, an isolating block~$B$, together
with its closed exit set~$B^-$, gives rise to an index pair, and the
relative homology~$H_*(B,B^-)$ is defined as the Conley index of the
largest invariant set~$S \subset B$. Conley theory implies that if
this homology is nontrivial, then~$S$ is necessarily nonempty. 
In other words, trivial flow through~$B$ is only possible for 
trivial Conley index. This distinction is mimicked by the 
multivector classification, since each multivector can be thought
of as a small isolating block, with the mouth taking the role
of the exit set.

Multivector fields turn out to be extremely flexible for encoding
dynamical behavior. This is partly due to the following result, 
which shows that multivector fields can be constructed uniquely
from a collection of {\em dynamical transitions\/}.

\begin{thm}
   Let $X$ be a Lefschetz complex and let $\mathcal{D}$ denote
    an arbitrary collection of subsets of $X$. Then there exists
    a uniquely determined minimal multivector field $\mathcal{V}$
    which satisfies the following:
    \begin{itemize}
    \item For every $D \in \mathcal{D}$ there is a
          $V \in \mathcal{V}$ with $D \subset V$.
    \end{itemize}
    Note that the sets in $\mathcal{D}$ do not have to be disjoint,
    and their union does not have to exhaust $X$. One can think of
    the sets in $\mathcal{D}$ as all allowable dynamical transitions.
\end{thm}

As we mentioned earlier, previous work, such as~\cite{boczko:etal:07a},
has relied on vector field transversality to drive the creation of phase
space decompositions suitable for dynamical analysis. A considerably
more flexible process involves the creation of a multivector field
based on the above result and an {\em arbitrary phase space decomposition\/}
via an underlying Lefschetz complex. This can always be accomplished in the
following way:
\begin{itemize}
   \item[(a)] Create a Lefschetz complex (e.g. simplicial, cubical) which
   discretizes the phase space region of interest.
   \item[(b)] For every cell $\sigma\in X$ with $\dim\sigma < \dim X$,
   use the vector field to find all higher-dimensional cells which can be
   reached directly from $\sigma$ to form sets $D_{\sigma} \subset X$.
   \item[(c)] Form the minimal multivector field $\cV$ associated with the
   collection $\mathcal{D} = \{ D_\sigma \, : \, \dim\sigma < \dim X \}$.
\end{itemize}
Thus, as long as one has an idea about the transitions that a system
has to be allowed to do, one can always find a smallest multivector field
realizing them. It is possible that the result leads to the trivial
multivector field $\mathcal{V} = \{ X \}$. In most cases, however,
the resulting multivector field is more useful, as will become clear 
from our examples in the last section. We emphasize that this method
does not require any transversality assumptions on the underlying
complex cells with respect to the flow of interest. In fact,
transversality is automatically created via the multivector field
construction, see~\cite{thorpe:wanner:p24a}. Therefore, it should 
be possible to use this construction using Lyapunov function
approximations.

%
%
\section{Combinatorial Dynamics and the Conley Index}

The previous section demonstrated how a multivector field can be 
constructed which encodes the continuous dynamics on the underlying
Lefschetz complex. In this section we examine the {\em combinatorial
dynamical system\/} induced by the multivector field which acts through
the iteration of a multivalued map. This {\em flow map\/} is given by
\begin{displaymath}
   \Pi_{\mathcal V}(x) = \mathrm{cl}\, x \cup [x]_{\mathcal V}
   \quad\text{ for all }\quad
   x \in X
\end{displaymath}
where $[x]_{\mathcal V}$ denotes the unique multivector in
${\mathcal V}$ which contains $x$. The definition of the
flow map shows that the induced dynamics combines two types
of behavior:
\begin{itemize}
   \item From a cell $x$, it is always possible to flow towards
   the boundary of the cell, i.e., to any one of its faces.
   \item In addition, it is always possible to move freely within
   a multivector.
\end{itemize}

In this way, the map encapsulates the possible dynamical 
behavior of the underlying classical dynamical system.

The multivalued map $\Pi_{\mathcal V} : X \multimap X$
naturally leads to a solution concept for multivector fields.
A {\em path\/} is a sequence $x_0, x_1, \ldots, x_n \in X$ such
that $x_k \in \Pi_{\mathcal{V}}(x_{k-1})$ for all indices
$k = 1,\ldots,n$.
A {\em solution\/} is then a map $\rho : \mathbb{Z} \to X$
which satisfies $\rho(k+1) \in \Pi_{\mathcal V}(\rho(k))$ for all
$k \in \mathbb{Z}$.
We say that this solution {\em passes through the cell\/}  $x \in X$
if in addition one has $\rho(0) = x$.
It is clear from the definition of the flow map that every
constant map is a solution, since we have the inclusion
$x \in \Pi_{\mathcal V}(x)$. Thus, it will be advantageous
to focus on a more restrictive notion of solution called an
{\em essential solution\/} where the underlying path has to leave a
regular multivector both in forward and in backward time.
\begin{defn}
   Let $\rho : \mathbb{Z} \to X$ be a solution for the
   multivector field $\mathcal{V}$. Then $\rho$ is
   an {\em essential solution\/}, if the following holds:

   If $\rho(k)$ lies in a regular multivector
   $V \in \mathcal{V}$ for $k \in \mathbb{Z}$, then there
   exist integers $\ell_1 < k < \ell_2$ with
   $\rho(\ell_i) \not\in V$ for $i = 1,2$.
\end{defn}

Next we introduce the notion of {\em invariance\/} in the context of
combinatorial dynamics. A subset $A\subset X$ of a Lefschetz complex
is {\em invariant\/} if, for every cell $x\in A$, there exists an
essential solution~$\rho$ through~$x$ which stays in~$A$. Conley
realized in \cite{conley:78a} that if one restricts the attention
to a more specialized notion of invariance, then topological methods
can be used to formulate a coherent general index theory for invariance.
Central to Conley's work is the concept of an {\em isolated invariant
set\/}, which has been adapted to the combinatorial setting.
\begin{defn}
   A closed set $N \subset X$ {\em isolates\/} an invariant set
    $S \subset N$, if the following two conditions are satisfied:
    \begin{itemize}
      \item Every path in $N$ with endpoints in $S$ is a path in
      $S$.
      \item We have $\Pi_{\mathcal{V}}(S) \subset N$.
    \end{itemize}
    An invariant set $S$ is an {\em isolated invariant set\/},
    if there exists a closed set $N$ which isolates $S$.
\end{defn}

One can show that if $N$ is an isolating set for an
isolated invariant set $S$, then any closed set
$S \subset M \subset N$ also isolates $S$. Thus, the closure
$\mathrm{cl}\, S$ is the smallest isolating set for $S$.
What follows is a very useful characterization due to \cite{lipinski:etal:23a}.
\begin{thm}
   An invariant set $S \subset X$ is an isolated invariant set,
   if and only if the following two conditions hold:
   \begin{itemize}
      \item $S$ is {\em $\mathcal{V}$-compatible\/}, i.e., it is the union
      of multivectors.
      \item $S$ is locally closed.
   \end{itemize}
   In this case, the isolated invariant set $S$ is isolated
   by its closure $\mathrm{cl}\, S$.
\end{thm}
Finally, since isolated invariant sets are locally closed,
we can define their {\em Conley index\/} as follows.
\begin{defn}
   Let $S \subset X$ be an isolated invariant set of the
   multivalued flow map $\Pi_{\mathcal{V}}$. Then the
   {\em Conley index of~$S$\/} is the relative (or Lefschetz)
   homology
   \begin{displaymath}
      CH_*(S) = H_*( \mathrm{cl}\, S, \mathrm{mo}\, S) \cong H_*(S)
   \end{displaymath}
\end{defn}

As in the classical case, the Conley index~$CH_*(S)$ encodes stability
information of the isolated invariant set~$S$. Moreover, the $k$-th Conley
index group~$CH_k(S)$ is in fact a vector space over~$F$, and therefore
determined completely by its dimension, i.e., the $k$-th Betti number.
   
%
%
\section{Morse Decompositions}

We now turn our attention to examining the global dynamics of the combinatorial
dynamical system through the notion of its {\em Morse decomposition\/} by extending
some concepts from classical dynamics.
\begin{defn}
   We say that an invariant set $S\subset X$ is an {\em attractor\/} if $\Pi_{\mathcal{V}}(S) = S$,
   and a {\em repeller\/} if $\Pi^{-1}_{\mathcal{V}}(S) = S$
\end{defn}
\begin{defn}
   The {\em $\mathcal{V}$-hull\/} of a set $A \subset X$ is the
   intersection of all $\mathcal{V}$-compatible and locally
   closed sets containing~$A$. It is denoted by
   $\langle A \rangle_{\mathcal{V}}$, and is the smallest 
   candidate for an isolated invariant set which contains $A$.
\end{defn}
\begin{defn}
   The $\alpha$- and $\omega$-limit sets of a solution~$\varphi$ are then defined as
   \begin{align}
      \alpha(\varphi) &=
      \left\langle \bigcap_{t \in \mathbb{Z}^-} \varphi\left( (-\infty,t] \right)   \right\rangle_{\mathcal{V}}
      \quad\text{ and } \nonumber \\
      \omega(\varphi) &=
      \left\langle \bigcap_{t \in \mathbb{Z}^+} \varphi\left( [t,+\infty) \right) \right\rangle_{\mathcal{V}}
      \nonumber
   \end{align}
\end{defn}

Note that the $\mathcal{V}$-hull of a set does not have to be invariant in general.
However, essential solutions present an important special case where both of their
limit sets are in fact isolated invariant sets.
\begin{thm}
   Let $\varphi$ be an essential solution in $X$. Then both
   limit sets $\alpha(\varphi)$ and $\omega(\varphi)$ are
   nonempty isolated invariant sets.
\end{thm}
The above notions allow us to decompose the global dynamics of
a multivector field. Loosely speaking, this is accomplished by
separating the dynamics into a recurrent part given by an indexed
collection of isolated invariant sets, and the gradient dynamics
between them. This can be abstracted through the concept of a
{\em Morse decomposition\/}.
    
\begin{defn}
   Assume that $X$ is an invariant set for the multivector
   field $\mathcal{V}$ and that $(\mathbb{P},\leq)$ is a
   finite poset. Then an indexed collection $\mathcal{M} =
   \left\{ M_p \, : \, p \in \mathbb{P} \right\}$ is called a
   {\em Morse decomposition\/} of $X$ if the following conditions are
   satisfied:
   \begin{itemize}
      \item The indexed family $\mathcal{M}$ is a family of mutually
      disjoint, isolated invariant subsets of $X$.
      \item For every essential solution~$\varphi$ we either have the
      inclusion
      $\mathrm{im} \, \varphi \subset M_r$ for an index $r \in \mathbb{P}$,
      or there exist two poset elements $p,q \in \mathbb{P}$ such
      that $q > p$, as well as $\alpha(\varphi) \subset M_q$
      and $\omega(\varphi) \subset M_p$.
   \end{itemize}
\end{defn}
The elements of $\mathcal{M}$ are called {\em Morse sets\/}.

There always exists a finest Morse decomposition $\mathcal{M}$ of a
given a combinatorial multivector field $\mathcal{V}$ on an
arbitrary Lefschetz complex $X$. It can be found by
determining the strongly connected components of the digraph
associated with the flow map $\Pi_{\mathcal{V}} : X \multimap X$
which contain essential solutions. See for example~\cite{lipinski:etal:23a}.

\begin{defn}
   The {\em connection set\/} $\mathcal{C}(A,B)$ from $A$ to $B$ is the
   set of cells $x \in X$ such that there exists an essential solution
   $\varphi$ through $x$ with $\alpha(\varphi) \subset A$ and $\omega(\varphi) \subset B$.
\end{defn}
The {\em connection set\/} captures the dynamics between two subsets $A,B \subset X$.
It is an isolated invariant set, and can be empty.
One can also combine connection set(s) with Morse set(s) to produce a special type of isolated
invariant sets for the multivector field called a {\em Morse interval\/}.
\begin{defn}
   Let $I \subset \mathbb{P}$ denote an interval in the index poset.
   Then the {\em Morse interval\/} is the collection
   \begin{equation}
      M_I \; = \; \bigcup_{p \in I} M_p \; \cup \;
                  \bigcup_{p,q \in I} \mathcal{C}( M_q, M_p )
   \end{equation}
\end{defn}

Morse intervals can be used to find attracting sets in a combinatorial
dynamical system. One can show that if the interval~$I$ is a {\em down set\/},
i.e., if with every $p \in I$ it also contains all elements $q \in \mathbb{P}$
with $q \le p$, then the Morse interval~$M_I$ is an attractor.

%
%
\section{Connection Matrices}

While a Morse decomposition represents the basic structure
of the global dynamics of a combinatorial dynamical system,
it does not directly provide more detailed information about
the dynamics between the Morse sets. For instance, we may be interested
in determining which of the associated connecting sets must be nonempty.
An algebraic way for making this determination is provided by the
{\em connection matrix\/} introduced in~\cite{franzosa:89a}.
The book by \cite{mrozek:wanner:25a} treats connection matrices
specifically in the setting of multivector fields and extends
the original definition.

Assume that we are given a Morse decomposition $\mathcal{M}$ of an
isolated invariant set $S$. Then the {\em connection matrix\/} is a
linear map defined on the direct sum of all Conley indices of the
Morse sets in the Morse decomposition.
\begin{displaymath}
   \Delta \; : \; \bigoplus_{q \in \mathbb{P}} CH_*(M_q)
   \to \bigoplus_{p \in \mathbb{P}} CH_*(M_p)
\end{displaymath}
One usually writes the connection matrix $\Delta$ as a matrix
in the form $\Delta = (\Delta(p,q))_{p,q \in \mathbb{P}}$,
which is indexed by the poset $\mathbb{P}$, and where
the entries $\Delta(p,q) : CH_*(M_q) \to CH_*(M_p)$ are
linear maps between the homological Conley indices, which are
vector spaces over~$F$. If~$I$ denotes an interval in the
poset~$\mathbb{P}$, then one further defines the restricted
connection matrix~$\Delta(I)$ as the minor~$(\Delta(p,q))_{p,q \in I}$.

Any connection matrix $\Delta$ has the following properties:
\begin{itemize}
\item $\Delta$ is {\em strictly upper triangular\/}, i.e., if
   $\Delta(p,q) \not= 0$ then one has $p < q$.
\item $\Delta$ is a {\em boundary operator\/}, i.e.,
   $\Delta^2 = 0$, and $\Delta$ maps
   $k$-th level homology to $(k-1)$-st level homology.
   \item For every interval $I$ in $\mathbb{P}$ we have
   \begin{displaymath}
      H_*\Delta(I) \; = \;
      \mathrm{ker}\, \Delta(I) / \mathrm{im}\, \Delta(I)
      \; \cong \; CH_*(M_I)
   \end{displaymath}
   In other words, the {\em Conley index of a Morse interval\/}
   can be determined via the {\em homology\/} of the associated
   {\em connection matrix minor\/} $\Delta(I)$.
   \item If $\{ p, q \}$ is an interval in $\mathbb{P}$ and
   $\Delta(p,q) \neq 0$, then the {\em connection set
   $\mathcal{C}(M_q,M_p)$ is not empty\/}.
\end{itemize}

We would like to point out that these properties do not
characterize connection matrices. In practice, a given
multivector field can have several different connection
matrices which can indicate the existence of bifurcations.
Uniqueness, however, can be observed if the underlying
system is a {\em gradient combinatorial Forman vector field\/} on
a Lefschetz complex. These are multivector fields in which
every multivector is either a singleton, and therefore
a critical cell, or a two-element Forman arrow. In
addition, a gradient combinatorial Forman vector field
cannot have any nontrivial periodic solutions, i.e., 
periodic solutions which are not constant and therefore
critical cells. See \cite{mrozek:wanner:25a}.

%
%
\section{Examples}

To illustrate the approach introduced in the last few sections
we close with some simple examples. These are meant to show the
potential of the approach, and are computed using the
Julia package {\tt ConleyDynamics.jl}, which can be found
at~\cite{conleydynamics}. This package is freely available
and contains implementations for all of the concepts described 
earlier. For more details on Julia we refer the reader
to~\cite{bezanson:17a}
%
%
\subsection{A Planar Gradient System}
\begin{figure}
   \begin{center}
      \includegraphics[width=7.0cm]{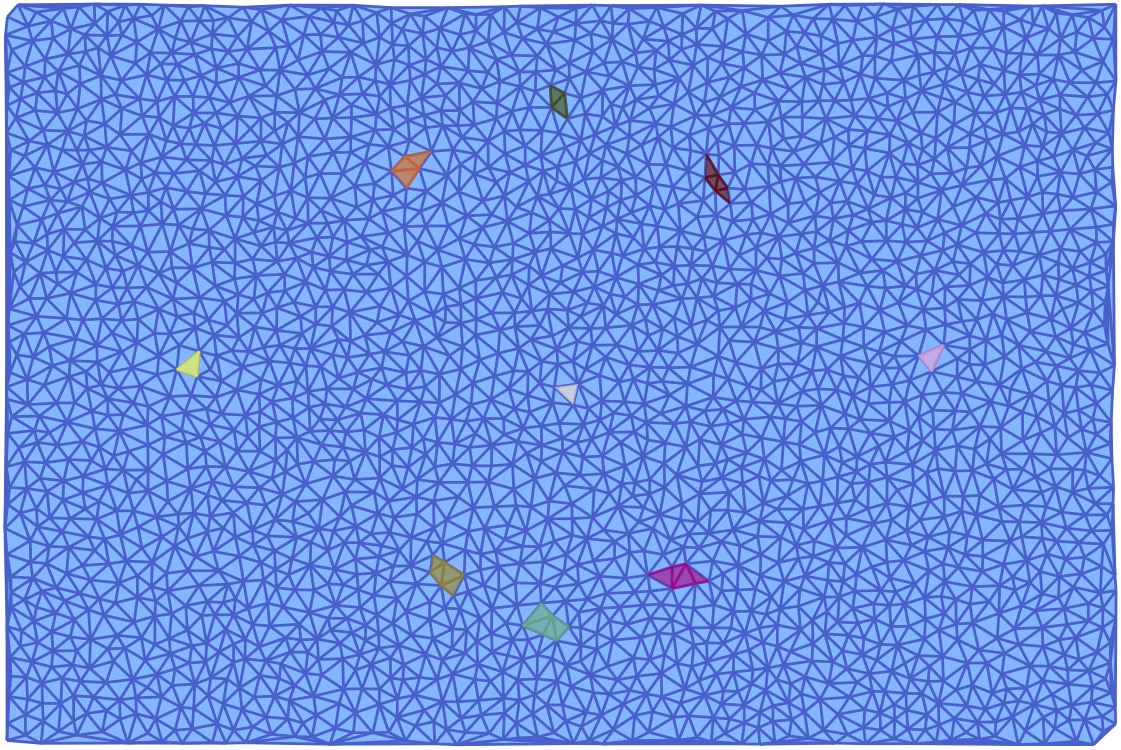}    
      \caption{Morse sets for the planar sample system, based
               on a random Delaunay triangulation.} 
      \label{fig:mvfanalysis1a}
   \end{center}
\end{figure}
\begin{figure}
   \begin{center}
      \includegraphics[width=7.0cm]{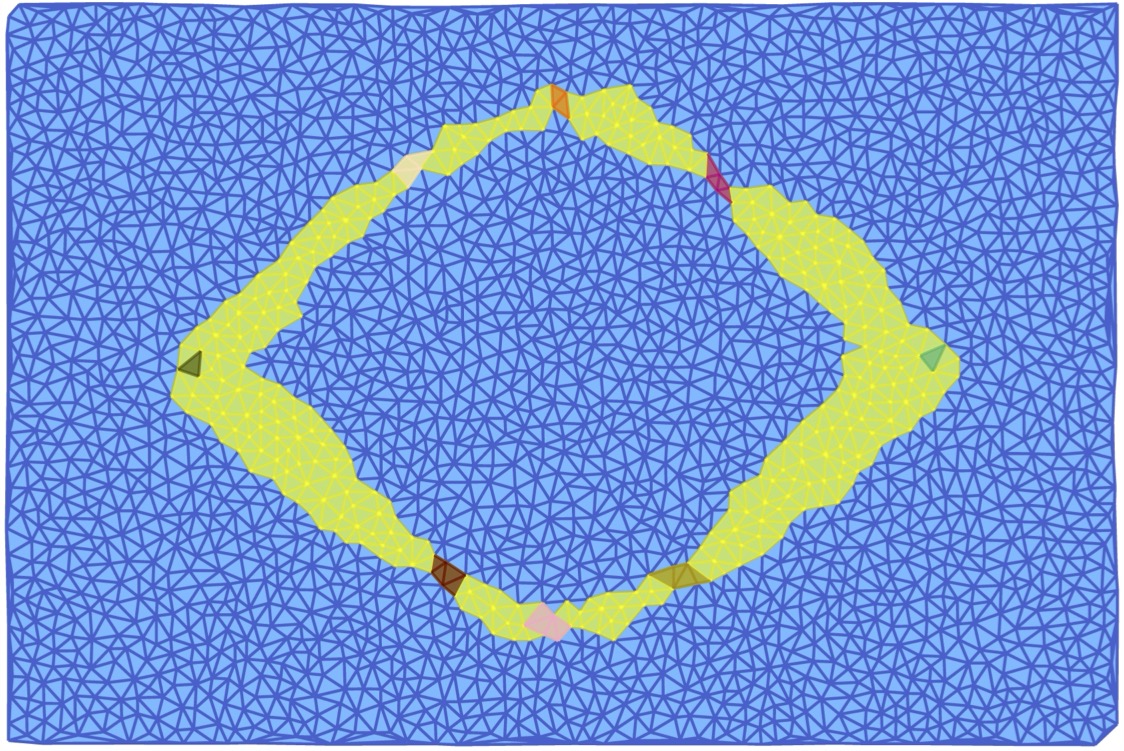}    
      \caption{Morse interval which includes all stable and
               index one stationary states of the planar sample
               system.} 
      \label{fig:mvfanalysis1b}
   \end{center}
\end{figure}
\begin{figure}
   \begin{center}
      \includegraphics[width=7.0cm]{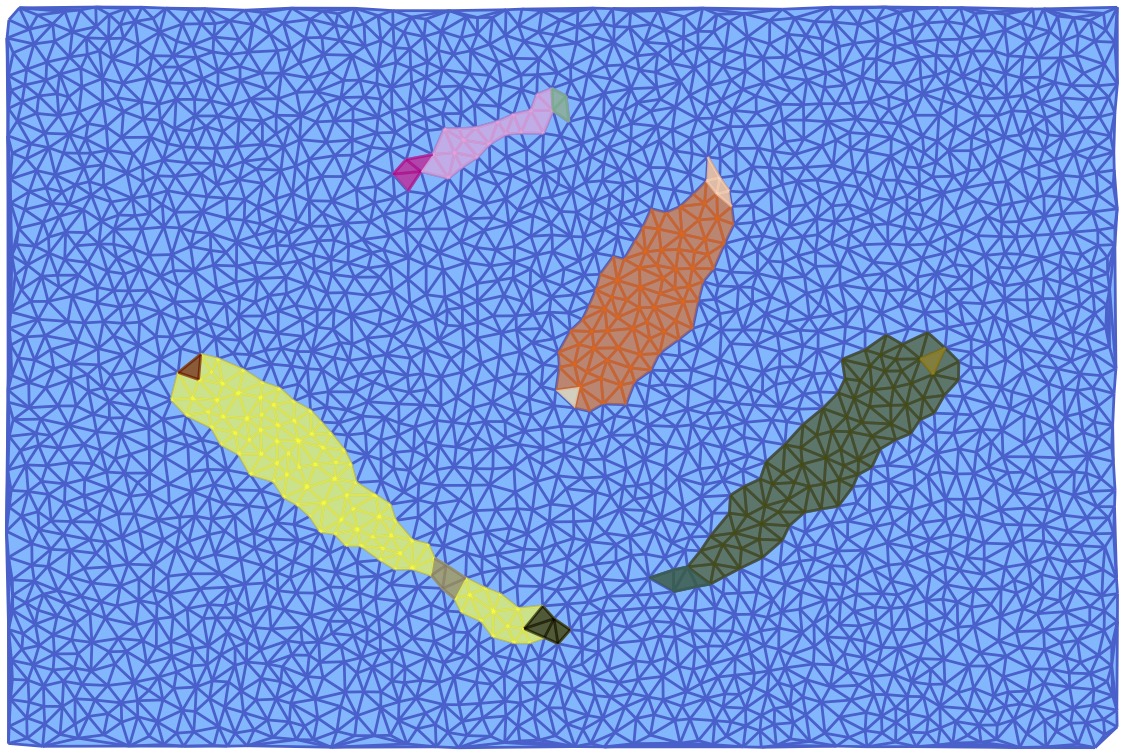}    
      \caption{Four sample Morse intervals for the planar
               system with nine equilibria.} 
      \label{fig:mvfanalysis1c}
   \end{center}
\end{figure}

Consider the gradient system in the plane defined by
\begin{align}
   \dot{x}_1 &= x_1 \left( 1 - x_1^2 - 3 x_2^2 \right) - 0.01
     \nonumber \\
   \dot{x}_2 &= x_2 \left( 1 - 3 x_1^2 - 2 x_2^2 \right) + 0.05
     \nonumber
\end{align}
This system has a global attractor, which is comprised of nine
equilibrium solutions and the connecting orbits between pairs
of them. There is an unstable equilibrium of index~$2$ close
to the origin, which is surrounded by four index~$1$ and four
asymptotically stable equilibria. The latter ones are close to
the coordinate axes.

Fig~\ref{fig:mvfanalysis1a} shows a coarse Delaunay triangulation
of a region containing the attractor. Based on the resulting Lefschetz
complex and the vector field, we then determine the induced multivector
field~$\mathcal{V}$ as described above. This field has a finest Morse
decomposition with nine isolated invariant sets, which are indicated 
by different colors in the figure, and each of which encloses one
equilibrium.

From the Morse decomposition, one can also determine Morse intervals.
The Morse intervals which contains all stable and index~$1$ equilibria
is shown in Fig~\ref{fig:mvfanalysis1b}, while four additional ones
are depicted in Fig~\ref{fig:mvfanalysis1c}. In all of these cases
the colored regions contain not only the respective stationary states
of the interval, but also enclosures for all connecting orbits between
them. In other words, the identified regions form rough outer
approximations for the heteroclinic orbits connecting these equilibria.

%
%
\subsection{The Three-Dimensional Allen-Cahn Projection}

As our second example we consider the three-dimensional system
\begin{align}
   \dot{x}_1 &= \left(\lambda-1\right)x_1 \nonumber \\ 
             &- \frac{3\lambda}{2\pi}\left[ x_1^3 - x_1^2 x_3 +
             x_2^2 x_3 + 2x_1 \left(x_2^2 + x_3^2 \right)\right]
             \nonumber \\
   \dot{x}_2 &= \left(\lambda-4\right)x_2 \nonumber \\
             &- \frac{3\lambda}{2\pi} x_2 \left[ 2x_1^2 + x_2^2 +
             2x_1 x_3 + 2x_3^2\right] \nonumber \\
   \dot{x}_3 &= \left(\lambda-9\right)x_3 \nonumber \\
             &- \frac{\lambda}{2\pi}\left[x_1\left(x_1^2 - 3x_2^2\right)
             - 3x_3 \left(2x_1^2 + 2x_2^2 + x_3^2\right)\right] \nonumber
\end{align}
This system is obtained from an Allen-Cahn type partial differential
equation via orthogonal projection onto the first three eigenfunctions
of the Laplacian. One can show that this system has seven equilibria
for the parameter value $\lambda = 3\pi$. One has index~$3$, and there
are pairs of equilibria for each of the indices~$2$, $1$, and~$0$.
\begin{figure}
   \begin{center}
      \includegraphics[width=7.0cm]{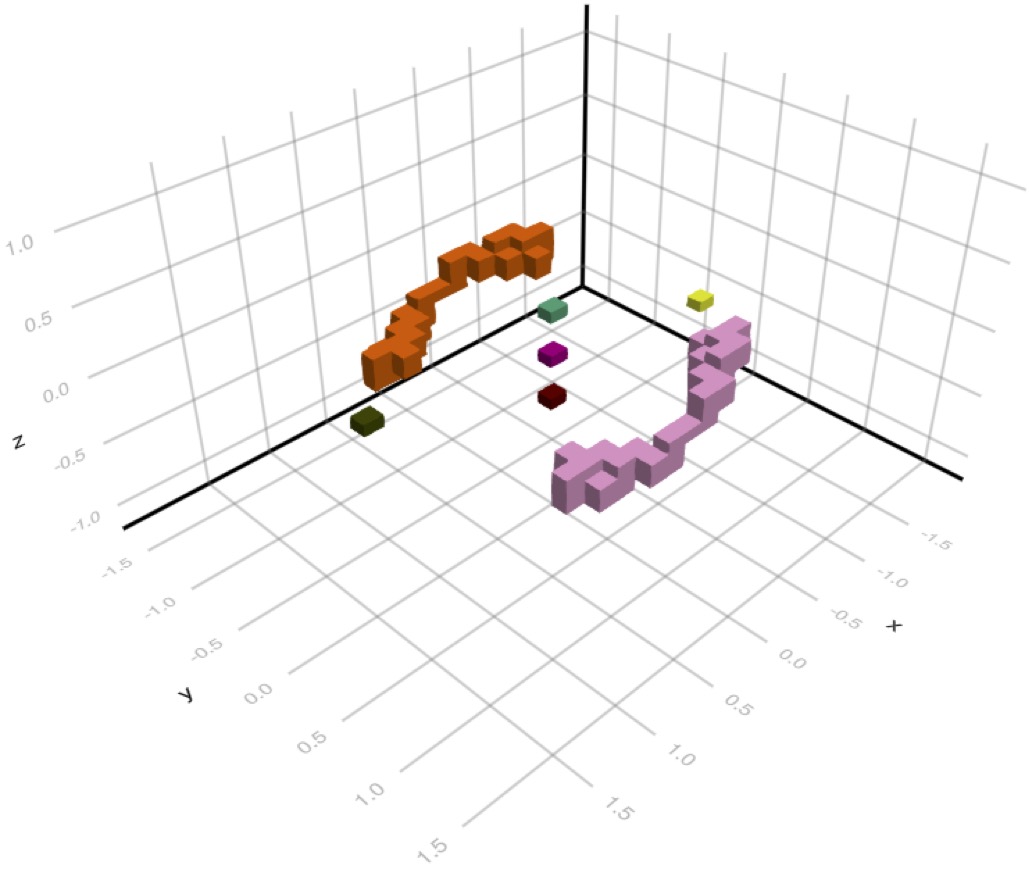}    
      \caption{Enclosures for the seven equilibrium solutions 
               of the three-dimensional sample system.} 
      \label{fig:allencahn3d_3_25}
   \end{center}
\end{figure}
\begin{figure}
   \begin{center}
      \includegraphics[width=7.0cm]{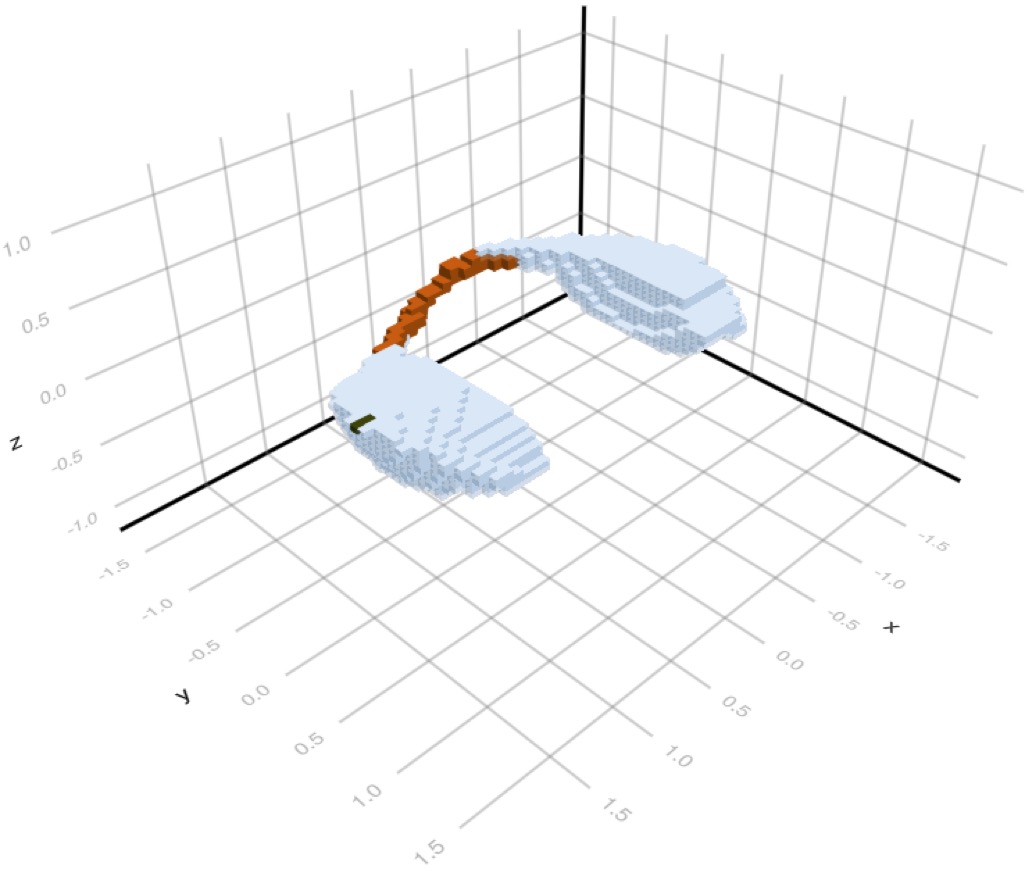}    
      \caption{A Morse interval for the three-dimensional sample
               system, which encloses a index~$1$ equilibrium, two stable
               stationary states, and the two heteroclinics between these
               states.} 
      \label{fig:acinterval3d_3_51a}
   \end{center}
\end{figure}

Fig~\ref{fig:allencahn3d_3_25} shows a coarse cubical decomposition of the
attracting region. The colored subsets are the isolated invariant sets
identified by the algorithm. The stationary states of index~$0$, $2$,
and~$3$ are all well-localized, but this cannot be said about the two
equilibria of index~$1$. The computed enclosures for the latter two are
elongated collections of cubes which are shown along the upper left and
lower right of the figure. This overestimation is a result of the use of
a strict cubical grid combined with the small discretization size.
Fig~\ref{fig:acinterval3d_3_51a} shows the view of the Morse interval
which corresponds to one of the index~$1$ equilibria and the two stable
stationary states. Note that this computation was performed with a
slightly finer discretization.

The above examples illustrate that the theory of multivector fields can be
used successfully to obtain statements regarding the classical dynamics of
an underlying system. By combining this approach with interval arithmetic
computations and connection matrix information, one can even derive
computer-assisted existence proofs for attractors and connecting orbits.
%
%
%
%
\bibliography{wanner1a,wanner1b,wanner2a,wanner2b,wanner2c}

\end{document}